\documentclass[12pt]{article}
\usepackage{amssymb}
\usepackage{amsmath,amsthm}
\usepackage{citesort}
\usepackage[latin1]{inputenc}

\usepackage{graphicx}
 \newtheorem{remark}{Remark}

 \newtheorem{lemma}[remark]{Lemma}
 \newtheorem{theorem}[remark]{Theorem}
 \newtheorem{proposition}[remark]{Proposition}
 \newtheorem{corollary}[remark]{Corollary}

\title{\vspace*{1.7cm} On the Randi\'{c} index and \\conditional parameters of a graph}

\author{J. A. Rodr\'{\i}guez\footnote{e-mail:\mbox{\tt
    juanalberto.rodriguez\@@uc3m.es}} and J. M. Sigarreta\footnote{e-mail:\mbox{\tt
    josemaria.sigarreta\@@uc3m.es}}\\{\em Department of Mathematics}\\ Carlos III of Madrid University\\
Avda. de la Universidad 30, 28911 Legan\'es (Madrid),  Spain}

\date{}

\begin{document}

\maketitle



\begin{center}
\end{center}

\begin{abstract}
The aim of this paper is to study some parameters of simple graphs
related with the degree of the vertices.   So, our main tool is
the $n\times n$ matrix ${\cal A}$ whose ($i,j$)-entry is
$$ a_{ij}= \left\lbrace \begin{array}{ll}
                         \frac{1}{\sqrt{\delta_i\delta_j}}  & {\rm if }\quad  v_i\sim v_j ;
                            \\ 0 &  {\rm otherwise,} \end{array}
                                                \right.  $$
where $\delta_i$ denotes the degree of the vertex $v_i$. We study
the Randi\'{c} index and some interesting particular cases of
conditional excess, conditional Wiener index, and conditional
diameter.  In particular, using the matrix ${\cal A}$ or its
eigenvalues, we obtain  tight bounds on the studied parameters.
\end{abstract}

\section{Introduction}

 In order to deduce properties of graphs from results and methods of algebra,
firstly we need to translate properties of graphs into algebraic
properties. In this sense,  a natural way is to consider algebraic
structures or algebraic objects as, for instance,  groups or
matrices. In particular, the use of matrices allows us to use
methods of linear algebra to derive properties of graphs. There
are various matrices that are naturally associated with graphs,
such as the adjacency matrix, the Laplacian matrix, and the
incidence matrix \cite{Biggs,Godsil,cvetkovic}. One of the main
aims of algebraic graph theory is to determine how, or whether,
properties of graphs are reflected in the algebraic properties of
such matrices \cite{Godsil}. The aim of this paper is to study the
Randi\'{c} index and some interesting particu\-lar cases of
conditional excess, conditional Wiener index, and conditional
diameter. All these parame\-ters are related with the degree of
the vertices of the graph. So, our main tool will be a suitable
adjacency matrix that we call \emph{degree-adjacency matrix}.

The plan of the paper is the following: in Section
\ref{degAdjMatrix}  we emphasize some of the main proper\-ties of
the degree-adjacency matrix. The remaining sections are devoted to
study the relationship between the degree-adjacency matrix (or its
eigenvalues) and several parameters of graphs. More precisely, in
Section \ref{randic} we obtain bounds on the Randi\'{c} index, in
Section \ref{ExcessAndWiener} we obtain bounds on a particular
case of conditional excess, Section \ref{DegreeDiam} is devoted to
bound the degree diameter and, finally, in section
\ref{CondWiener} we obtain bounds on a particular case of
conditional Wiener index.

We begin by stating some notation. In this paper all graphs
$\Gamma=(V,E)$ will be finite, undirected and simple. We will
assume that $|V|=n$ and $|E|=m$. The distance between  vertices
$u, v\in V(\Gamma)$ will be denoted by $\partial(u,v)$. The degree
of a vertex $v_i\in V(\Gamma)$ will be denoted by $\delta(v_i)$
(or by $\delta_i$ for short), the minimum degree of $\Gamma$ will
be denoted by $\delta$ and the maximum by $\Delta$.

\section{Degree-adjacency matrix}\label{degAdjMatrix}

 We define the $\emph{degree-adjacency matrix}$ of a
graph $\Gamma$ of order $n$ as the $n\times n$ matrix ${\cal A}$
whose ($i,j$)-entry is
\[ a_{ij}= \left\lbrace \begin{array}{ll} \frac{1}{\sqrt{\delta_i\delta_j}}  & {\rm if }\quad  v_i\sim v_j ;
                                                    \\
                                                0 &  {\rm otherwise.} \end{array}
                                                \right. \]

The matrix ${\cal A}$ can be regarded as the adjacency matrix of a
weighted graph in which the edge-weight ${\emph{R}}(v_iv_j)$ of
the edge $v_iv_j$ is equal to
${\emph{R}}(v_iv_j)=\frac{1}{\sqrt{\delta_i\delta_j}}$, thus
justifying the terminology used. The weight ${\emph{R}}(v_iv_j)$
will be called the {\em Randi\'{c} weight} of the edge $v_iv_j\in
E $. We will say that a graph is {\em weight-regular} if each of
its edges has the same Randi\'{c} weight. Particular cases of
weight-regular graphs are the class of regular graphs and the
class of  semi-regular bipartite graphs.

If we consider the vector $\nu=(\sqrt{\delta_1}, \sqrt{\delta_2},
..., \sqrt{\delta_n})$, then we have ${\cal A} \nu =\nu$. Thus,
$\lambda=1$ is an eigenvalue of ${\cal A}$ and $\nu$ is an
eigenvector associated to $\lambda$. Hence, as ${\cal A}$ is
non-negative and irreducible in the case of connected graphs, by
the Perron-Frobenius theorem,  $\lambda=1$ is a simple eigenvalue
and $\lambda=1\ge | \lambda_j |$ for every eigenvalue $\lambda_j$
of ${\cal A}$. Therefore, we have

\begin{equation}\label{contraida}
\|{\cal A}x\|\le \|x\|, \quad \forall x\in \mathbb{R}^n.
\end{equation}
Notice that the above inequality holds also in the case of
non-connected graphs.

Hereafter the eigenvalues of ${\cal A}$ will be called {\em
degree-adjacency eigenvalues of} $\Gamma$.

\begin{figure}[h]
\begin{center}
\caption{Two cospectral graphs but not cospectral with regard to
${\cal A}$} \label{ej1}
\includegraphics[width=0.3\textwidth]{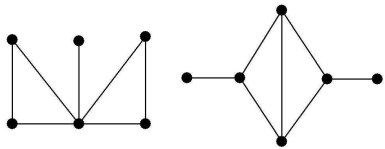}
\end{center}
\end{figure}

It is well-known that there are non-isomorphic graphs that have
the same standard adjacency eigenvalues with the same
multiplicities (the so called cospectral graphs). For instance,
two connected graphs, both having the characteristic polynomial
$P(x)=x^6-7x^4-4x^3+7x^2+4x-1$, are shown in Figure \ref{ej1}.
\vspace{-0,5cm} Therefore, we can try to study cospectral graphs
by using an alternative matrix, for instance, the degree-adjacency
matrix ${\cal A}$. If we consider the matrix ${\cal A}$, the
eigenvalues of both graphs are different: the left hand side graph
has degree-adjacency eigenvalues 1, $\pm\frac{1}{2}$ and
$-\frac{1}{4}\left(1\pm \sqrt{2.6} \right)$ (where the eigenvalue
$-\frac{1}{2}$ has multiplicity 2), on the other hand, the right
hand side graph has degree-adjacency eigenvalues 1, $\frac{-1\pm
\sqrt{2}}{3}$, $\pm \frac{\sqrt{3}}{3}$ and $- \frac{1}{3}$. Even
so, the degree-adjacency eigenvalues do not determine the graph.
That is, there are non-isomorphic graphs (and non-cospectral) that
are cospectral with regard to the degree-adjacency matrix. For
instance, the degree-adjacency eigenvalues of the cycle graph
$C_4$ and the  semi-regular bipartite graph $K_{1,3}$ are the
same: $1,0,0,-1$. However, the standard eigenvalues  are
$2,0,0,-2$, in the case of $C_4$, and $\sqrt{3},0,0,-\sqrt{3}$ in
the case of $K_{1,3}$.

 It is easy to see that there are some classes
of graphs in which the standard eigenvalues,
$\vartheta_1\ge\vartheta_2\ge\cdots\ge\vartheta_n$, and the
degree-adjacency eigenvalues,
$\lambda_1\ge\lambda_2\ge\cdots\ge\lambda_n$, are directly
related. For instance, in the case of weight-regular graphs, of
weight $w^{-1}$, the adjacency matrix, ${\bf A}$, and the
degree-adjacency matrix are related by ${\cal
A}(\Gamma)=\frac{1}{w}{\bf A}(\Gamma)$. Thus, the eigenvalues of
both matrices are related by
\begin{equation}
\lambda_l=\frac{\vartheta_l}{w}, \quad l\in\{1,...,n\}.
\end{equation}

As in the case of the adjacency matrix, there are some classes of
graphs in which we can deduce a formula to compute the
characteristic polynomial, $\Psi$, of the degree-adjacency matrix.
For instance, from the degree-adjacency matrix of the path graph,
$\Gamma=P_n$, we deduce that
$$\Psi({\cal A}(P_n),\lambda)=\lambda \Phi_{n}(\lambda)+\frac{1}{2}\Phi_{n-1}(\lambda), \quad n\ge 3,$$
where $$\Phi_{n}(\lambda)=-\lambda
\Phi_{n-1}(\lambda)-\frac{1}{4}\Phi_{n-2}(\lambda),\quad
\Phi_2(\lambda)=-\lambda \quad {\rm and} \quad
\Phi_3(\lambda)=\lambda^2-\frac{1}{2}. $$

Hereafter,  in the general case of an arbitrary graph, we will
consider that the characteristic polynomial, $\Psi({\cal A},
\lambda)=det(\lambda {\bf I}-{\cal A})$, is of the form
$$\Psi({\cal A}, \lambda)=\lambda^n+c_1\lambda^{n-1}+c_2\lambda^{n-2}+\cdots +c_{n-1}\lambda+c_n.$$
We can compute the first coefficients of $\Psi$ by using a
well-known result of theory of matrices: all the coefficients can
be expressed in terms of the principal minors of ${\cal A}$.


\begin{proposition}
Let $\Gamma$ be a graph. The coefficients of the characteristic
polynomial of ${\cal A}={\cal A}(\Gamma)$ satisfy:
\begin{equation}
c_1=0;
\end{equation}
\begin{equation} \label{c2}
-c_2=\frac{1}{2}\sum_{l=1}^{n}\lambda_l^2=\sum_{i\sim
j}\frac{1}{\delta_i\delta_j};
\end{equation}

\begin{equation}\label{c3}
c_3=\sum_{\langle i,j,k \rangle \simeq
K_3}\frac{-2}{\delta_i\delta_j\delta_k},
\end{equation}
where $\langle i,j,k \rangle \simeq K_3$ runs over all subgraphs
of $\Gamma$ induced by  $\{v_i,v_j,v_k\}$  and isomorphic to
$K_3$.
\end{proposition}

\begin{proof}
For each $r=1,2,...,n$, the number $(-1)^rc_r$ is the sum of those
principal minors of ${\cal A}$ which have order $r$. Thus, we
derive the result as follows. Since ${\cal A}$ has diagonal
entries all zero, $c_1=0.$  A principal non-null minor of order 2
must be of the form
$$\left|
 \begin{array}{cc}
0 & \frac{1}{\sqrt{\delta_i\delta_j}}  \\
\frac{1}{\sqrt{\delta_i\delta_j}}& 0  \\
\end{array}
\right|=\frac{-1}{\delta_i\delta_j} .$$ There is one such minor
for each edge of $\Gamma$. Moreover, since the trace of a square
matrix is also equal to the sum of its eigenvalues, we have
$$\sum_{l=1}^{n}\lambda_l^2=tr \left({\cal A}^2\right)=2\sum_{i\sim
j}\frac{1}{\delta_i\delta_j}.$$ Thus, (\ref{c2}) follows. On the
other hand, the only non-null principal minor of order 3 is
$$\left|
 \begin{array}{ccc}
0 & \frac{1}{\sqrt{\delta_i\delta_j}} & \frac{1}{\sqrt{\delta_i\delta_k}}\\
\frac{1}{\sqrt{\delta_j\delta_i}}& 0  &\frac{1}{\sqrt{\delta_j\delta_k}}\\
\frac{1}{\sqrt{\delta_k\delta_i}}&\frac{1}{\sqrt{\delta_k\delta_j}}
&0
\end{array}
\right|=\frac{2}{\delta_i\delta_j\delta_k} .$$ There is one such
minor for each triangle of $\Gamma$. Hence, (\ref{c3}) follows.
\end{proof}

Notice that the coefficient $c_2$ is immediately bounded from
(\ref{c2}):
\begin{equation}
\frac{m}{\Delta^2}\le -c_2\le \frac{m}{\delta^2}.
\end{equation}

\begin{corollary}
A graph is regular if, and only if, its order is $2|c_2|\Delta$.
\end{corollary}


In Section \ref{randic} we will show the relationship between
$c_2$ and the generalized Randi\'{c} index.

We remark that  the spectrum of ${\cal A}$ can be computed
directly from the adjacency matrix ${\bf A}$ and the degree
sequence. That is,
\begin{equation}
det({\cal A}-\lambda {\bf I})\cdot \prod_{j=1}^n \delta_j=det({\bf
A}-\lambda {\bf D}),
\end{equation}
where ${\bf D}=diag(\delta_1,\delta_2,...,\delta_n)$ is the
diagonal matrix whose diagonal entries are the degrees of the
vertices of $\Gamma$.

 There are other properties of the degree-adjacency
matrix that have been obtained previously (see, for instance
\cite{cvetkovic}),  in the following theorem we cite some of them.

\begin{theorem} {\rm  \cite{cvetkovic}}. \mbox{ }

\begin{itemize}
\item The number of connected components of $\Gamma$ is equal to
the multiplicity of the eigenvalue $1$ of ${\cal A}$.

\item Let $\Gamma$ be a graph without isolated vertices. $\Gamma$
is bipartite if and only if
$\Psi(\Gamma,\lambda)=\Psi(\Gamma,-\lambda)$.

\item Let $\Gamma$ be a connected graph. $\Gamma$ is bipartite if
and only if $-1$ is an eigenvalue of ${\cal A}$.
\end{itemize}
\end{theorem}

 We identify the degree-adjacency matrix
${\cal A}$ with an endomorphism of the ``vertex-space" of
$\Gamma$, $l^2(V(\Gamma))$ which, for any given indexing of the
vertices, is isomorphic to $\mathbb{R}^n$. Thus, for any vertex
$v_i \in V(\Gamma)$, $e_i$ will denote the corresponding  unit
vector of the canonical base of $\mathbb{R}^n$.

If for two vertices $v_i,v_j\in V(\Gamma)$ we have $\partial
(v_i,v_j)>k$ then $({\cal A}^k(\Gamma))_{ij}=0$. Thus, for a real
polynomial $P$ of degree $k$, we have
\begin{equation}
\partial (v_i,v_j)>k \Rightarrow  P({\cal A}(\Gamma))_{ij}=0 .\label{just}
\end{equation}
Through this fact  we will study some metric parameters of graphs.

\section{Randi\'{c} index} \label{randic}

The $\emph{Randi\'{c} index}$, $\emph{R}(\Gamma)$, of a graph
$\Gamma$ was introduced by the chemist Milan Randi\'{c} in 1975
\cite{Randic} as
\[
{\emph{R}}(\Gamma)=\sum_{v_i\sim v_j}
\frac{1}{\sqrt{\delta_i\delta_j}}.
\]
 This topological index, sometimes called ${\emph{connectivity
index,}}$ has been success\-fully related to physical and chemical
properties of organic molecules and became one of the most popular
molecular descriptors.

The Randi\'{c} index has the following trivial bounds:
\begin{equation}
\frac{m}{\Delta}\le {\emph{R}}(\Gamma)\le \frac{m}{\delta}.
\end{equation}
 Equality holds if, and only if, $\Gamma$ is regular. Moreover, there are non-trivial bounds as the following
\cite{Bollobas}:
\begin{equation}\label{BollErdos}
\sqrt{n-1}\le {\emph{R}}(\Gamma)\le \frac{n}{2}.
\end{equation}
Equality on the right-hand side holds if, and only if, $\Gamma$ is
a graph whose all components are regular of (not necessarily
equal) degrees greater than zero. Equality on the left-hand side
holds if, and only if, $\Gamma$ is a star \cite{Bollobas}.

We emphasize that the degree-adjacency matrix allows us to obtain
a short proof of the right hand side of (\ref{BollErdos}):  by the
Cauchy-Schwarz inequality and (\ref{contraida}) we have
 $$2{\emph{R}}(\Gamma)=\langle {\cal A}{\bf j},{\bf
j}\rangle\le \|{\cal A}{\bf j}\|\|{\bf j}\|\le\|{\bf j}\|^2=n,$$
where ${\bf j}=(1,1,...,1)\in \mathbb{R}^n$.

The {\em zeroth-order} Randi\'{c} index  is defined as
\[
{\emph{R}}_0(\Gamma)=\sum_{v\in V(\Gamma)}
\frac{1}{\sqrt{\delta(v)}}.
\]
Trivially, ${\emph{R}}_0(\Gamma)$ is bounded by
\begin{equation}
\frac{n}{\sqrt{\Delta}}\le {\emph{R}}_0(\Gamma)\le
\frac{n}{\sqrt{\delta}}.
\end{equation}
The equality holds if, and only if, $\Gamma$ is regular of degree
greater than zero.

The Randi\'{c} index has been generalized \cite{gutman} as
\[
{\emph{R}}_\alpha(\Gamma)=\sum_{v_i\sim v_j}
\left(\delta_i\delta_j\right)^\alpha , \quad  \alpha\ne 0.
\]
Obviously, the standard Randi\'{c} index is obtained when
$\alpha=-\frac{1}{2}$.

In the chemical literature the quantity
$${\emph{R}}_1(\Gamma)=\sum_{v_i\sim v_j}\delta_i\delta_j$$
is called  {\em the second Zagreb index} \cite{Zagreb}. The second
Zagreb index was  bounded in \cite{Bollobas} by
\begin{equation}
{\emph{R}}_1(\Gamma)\le m \left(\frac{\sqrt{8m+1}-1}{2}\right)^2.
\end{equation}
Moreover, by (\ref{c2}) we have
\begin{equation}\label{c22}
{\emph{R}}_{-1}(\Gamma)=|c_2|.
\end{equation}

The $\emph{higher-order Randi\'{c} index}$ or $\emph{higher-order
connectivity index}$ is  also of  interest in mole\-cular graph
theory. For $t\ge 1$, the higher-order Randi\'{c} index is defined
as
\[
{\emph{R}}^{(t)}(\Gamma)=\sum_{v_{i_{1}}- v_{i_{2}}-\cdots -
v_{i_{t+1}}} \frac{1}{\sqrt{\delta_{i_1}\delta_{i_2}\cdots
\delta_{i_{t+1}}}},
\]
where $v_{i_{1}}- v_{i_{2}}-\cdots - v_{i_{t+1}}$ runs over all
paths of length $t$ in $\Gamma$.

Now we are going to obtain tight bounds on
${\emph{R}}_\alpha(\Gamma)$. Moreover, we are going to obtain
tight bounds on ${\emph{R}}_\alpha(\Gamma)$, $\alpha\neq -1, 0$,
and ${\emph{R}}^{(2)}(\Gamma)$ in terms of $c_2$ (the coefficient
of $\lambda^{n-2}$ in the characteristic polynomial,$\Psi({\cal
A}, \lambda)$, of the degree-adjacency matrix of $\Gamma$).

\begin{theorem} Let $\Gamma$ be a simple graph of
order $n$ and size  $m$.

\begin{itemize}
\item[(a)]{ The zeroth-order Randi\'{c} index is bounded by
$$\frac{n^3}{2m}\le {\emph{R}}_0^2(\Gamma).$$
The equality holds if, and only if, $\Gamma$ is regular.}
\item[(b)]{Let $\alpha_1,\alpha_2\in \mathbb{R}\setminus\{0\}$
such that $\alpha_1<\alpha_2$. Then
\begin{equation}\label{cotsup}
\alpha_1\alpha_2>0\Rightarrow
{\emph{R}}_{\alpha_1}^{\alpha_2}(\Gamma)m^{\alpha_1}\le
{\emph{R}}_{\alpha_2}^{\alpha_1}(\Gamma)m^{\alpha_2}.
\end{equation}
\begin{equation}\label{cotinf}
\alpha_1\alpha_2<0\Rightarrow
{\emph{R}}_{\alpha_2}^{\alpha_1}(\Gamma)m^{\alpha_2}\le
{\emph{R}}_{\alpha_1}^{\alpha_2}(\Gamma)m^{\alpha_1}.
\end{equation}
 The equalities hold
if, and only if, $\Gamma$ is weight-regular.}

 \item[(c)]{Let
$\vartheta_1\ge\vartheta_2\ge\cdots\ge\vartheta_n$ be the standard
eigenvalues of $\Gamma$ and let
$\lambda_1\ge\lambda_2\ge\cdots\ge\lambda_n$ be the
degree-adjacency eigenvalues of $\Gamma$, then
$${\emph{R}}(\Gamma)\le
\frac{1}{2}\sum_{i=1}^n|\lambda_i\vartheta_i|.$$}
\end{itemize}
\end{theorem}

\begin{proof}

\mbox{ }

\begin{itemize}
\item[(a)]{Application of   the Jensen's inequality to the convex
function $f(x)=x^{-2}$ leads to the result. That is
$$\frac{n^2}{{\emph{R}}_0^2(\Gamma)}=f\left(\frac{{\emph{R}}_0(\Gamma)}{n}\right)\le
 \displaystyle\frac{1}{n}\sum_{v_i\in
 V(\Gamma)}\delta_i=\frac{2m}{n}.$$}

 \item[(b)]{Let $g(x)=x^{\frac{\alpha_2}{\alpha_1}}$, where  $x>0$. If
 ($\alpha_1<0$ and $\alpha_2>0$) or ($0<\alpha_1<\alpha_2$),
application of the Jensen's inequality to the convex function $g$
leads to
\begin{equation}\label{jensen}
\left(\frac{{\emph{R}}_{\alpha_1}(\Gamma)}{m}\right)^{\frac{\alpha_2}{\alpha_1}}
=
g\left(\frac{{\emph{R}}_{\alpha_1}(\Gamma)}{m}\right)\le\frac{{\emph{R}}_{\alpha_2}(\Gamma)}{m}.
\end{equation}
Thus, by (\ref{jensen}), if $\alpha_1<0$ and $\alpha_2>0$ we
obtain
\begin{equation}\label{cota1}
{\emph{R}}_{\alpha_2}^{\alpha_1}(\Gamma)m^{\alpha_2}\le
{\emph{R}}_{\alpha_1}^{\alpha_2}(\Gamma)m^{\alpha_1}
\end{equation}
and, if $0<\alpha_1<\alpha_2$, we obtain
\begin{equation} \label{cota2}
{\emph{R}}_{\alpha_1}^{\alpha_2}(\Gamma)m^{\alpha_1}\le
{\emph{R}}_{\alpha_2}^{\alpha_1}(\Gamma)m^{\alpha_2}.
\end{equation}
Analogously, if $\alpha_1<\alpha_2<0$, application of the Jensen's
inequality to the concave function $g$ leads to (\ref{cota2}).
Hence, the result follows.}

 \item[(c)]{The result is obtained by $2{\emph{R}}(\Gamma)=Tr({\cal A}{\bf
A})\le \displaystyle\sum_{i=1}^n|\lambda_i\vartheta_i|.$}
 \end{itemize}
\end{proof}

Notice that, in the case of weight-regular graphs, the  bound (c)
is attained. Moreover, as a particular case of (b), by
(\ref{c22}), we deduce de following result.

\begin{corollary}
Let $\Gamma$ be a simple graph of size  $m$. Then
\begin{equation}\label{cotinf}
\alpha \in \mathbb{R}\setminus[-1,0] \Rightarrow
\frac{m^{\alpha+1}}{|c_2|^\alpha}\le {\emph{R}}_\alpha(\Gamma)
\end{equation}
\begin{equation}\label{cotsup}
\alpha \in (-1,0)\Rightarrow {\emph{R}}_\alpha(\Gamma)\le
\frac{m^{\alpha+1}}{|c_2|^\alpha}
\end{equation}
 The equalities hold
if, and only if, $\Gamma$ is weight-regular.
\end{corollary}

As a particular case of above corollary we obtain
\begin{equation}
{\emph{R}}(\Gamma)\le \sqrt{m|c_2|}.
\end{equation}

\begin{theorem} \label{Th1}
Let $\Gamma$ be a simple and connected graph of order $n$ and size
$m$. Let $\phi$ denotes de graph invariant defined as
$\phi=\left(\sum_{i=1}^n \sqrt{\delta_i} \right)^2 / 2m$, and let
$\delta_*=\displaystyle\min_{\delta_j>1}\{\delta_j\}$. Then
\[ \left(\frac{(2R(\Gamma)-\phi)^2}{2(n-\phi)} +\frac{\phi}{2}+c_2 \right)\sqrt{\delta_*}\le
R^{(2)}(\Gamma)\le \sqrt{\Delta}\left(\frac{n}{2}+c_2\right),\]
 where the lower bound holds only in the case of a non-regular
 graph.
\end{theorem}

\begin{proof}
For the  vector ${\bf j}=(1,1,...,1)\in \mathbb{R}^n$ we consider
the following decomposition
\begin{equation}\label{decomposition}
{\bf j}=\frac{\langle {\bf j},\nu\rangle}{\|\nu
\|^2}\nu+z=\frac{\sum_{i=1}^n \sqrt{\delta_i} }{ \sum_{i=1}^n
\delta_i }\nu+z,
\end{equation}
 where $z\in\nu^\perp . $
 Then we have
\begin{align*}
 2\emph{R}(\Gamma)&=\langle {\cal A}{\bf j},{\bf j}\rangle  \\
                    &=\left\langle \frac{\sum_{i=1}^n \sqrt{\delta_i}  }{
\sum_{i=1}^n \delta_i }\nu,\frac{\sum_{i=1}^n \sqrt{\delta_i}  }{
\sum_{i=1}^n \delta_i }\nu \right\rangle + \left\langle {\cal
A}z,z \right\rangle \\
 &=\frac{\left(\sum_{i=1}^n \sqrt{\delta_i} \right)^2 }{
\sum_{i=1}^n \delta_i }+\left\langle {\cal A}z,z \right\rangle \\
&=\phi+\left\langle {\cal A}z,z \right\rangle .
\end{align*}
Thus, $2\emph{R}(\Gamma)-\phi= \left\langle {\cal A}z,z
\right\rangle$ and by the Cauchy-Schwarz inequality we obtain $
|2\emph{R}(\Gamma) -\phi|\le \| {\cal A}z \| \| z\| $ and from $
\|z\|=\sqrt{n-\phi}$ and $ \|{\cal A}z \|=\sqrt{\|{\cal A}{\bf j}
\|^2-\phi}$ we obtain
\begin{equation} \label{ineq}
|2\emph{R}(\Gamma) - \phi|\le\sqrt{\left(\|{\cal A}{\bf j}
\|^2-\phi\right)(n-\phi)}.
\end{equation}
Moreover,
\begin{align}
\|{\cal A}{\bf j} \|^2&=2\sum_{v_i\sim
v_j}\frac{1}{\delta_i\delta_j}+2\sum_{v_i-v_j-v_k}\frac{1}{\sqrt{\delta_j\delta_i\delta_j\delta_k}}\label{Aj}\\
&\le 2\sum_{v_i\sim
v_j}\frac{1}{\delta_i\delta_j}+\frac{2}{\sqrt{\delta_*}}\sum_{v_i-v_j-v_k}\frac{1}{\sqrt{\delta_i\delta_j\delta_k}}.
\end{align}
Hence, by (\ref{c2}) we obtain
$$\|{\cal A}{\bf j} \|^2\le 2\left(\frac{\emph{R}^{(2)}(\Gamma)}{\sqrt{\delta_*}}-c_2\right).$$
Thus, if $\Gamma$ is non-regular, by the above inequality and
(\ref{ineq}) we conclude the proof of the left hand side
inequality. On the other hand, by (\ref{contraida}) we have
$\|{\cal A}{\bf j} \|^2\le \|{\bf j}\|^2$=n, then, by (\ref{Aj})
and (\ref{c2}) we have
$$-2c_2+\frac{2\emph{R}^{(2)}(\Gamma)}{\sqrt{\Delta}}  \le -2c_2+2\sum_{v_i-v_j-v_k}\frac{1}{\sqrt{\delta_j\delta_i\delta_j\delta_k}}\le
n. $$ Hence, the result follows.
\end{proof}

The above  bounds are attained, for instance, in the case of the
star graphs. Moreover, the upper bound is attained also in the
case of regular graphs. The reader is referred to \cite{Randic1}
for a complementary study on the Randi\'{c} index.

\section{Conditional excess}
\label{ExcessAndWiener}

Let $D(\Gamma)$ denotes the diameter of $\Gamma$. We define, for
any $k=0,1,\dots, D(\Gamma)$, the {\em $k$-excess} of a vertex
$u\in V(\Gamma)$, denoted by ${\bf e}_k(u)$, as the number of
vertices which are at distance greater than $k$ from $u$. That is,
$$ {\bf e}_k(u) = \vert \lbrace v\in V: \partial (u,v)>k\rbrace \vert.$$
Then, trivially, ${\bf e}_0(u)=n-1$,  ${\bf e}_{D(\Gamma)}(u) =
{\bf e}_{\varepsilon(u)}(u) = 0$ and ${\bf e}_k(u) = 0$
 if and only if $\varepsilon(u)\le k $,  where  $\varepsilon(u)$  denotes the eccentricity of $u$.
 The name ``{\em excess}" is borrowed from Biggs  \cite{Biggs},
 in which he gives a lower bound, in terms of the adjacency eigenvalues of a graph, for the excess
 ${\bf e}_r(u)$ of any vertex $u$ in a $\delta$-regular graph with odd girth $g=2r+1$.
 The excess of a vertex was studied by Fiol and Garriga    \cite{from-local}
 using the adjacency eigenvalues of a graph, and by Yebra and the first author of this paper in \cite{eama}
  using the Laplacian eigenvalues.

The $k$-{\em excess} of $\Gamma$, denoted by  ${\bf e}_k$, is
defined as
$$ {\bf e}_k = \max_{v_i\in V(\Gamma)}\lbrace {\bf e}_k(v_i)\rbrace.$$
This parameter  was studied by Yebra and the first author of this
paper  in \cite{exceso} using the Laplacian spectrum and the
$k$-alternating polynomials.

We define the {\em conditional excess} of a vertex $v\in
V(\Gamma)$ as follows:
$${\bf e}_k^\wp(u) :=\vert \lbrace v\in \wp: \partial (u,v)>k  \rbrace\vert ,$$
where $\wp$ is a property of some vertices of $\Gamma$ and $v\in
\wp$ means that the vertex  $v$ satisfies the property $\wp$. In
this section we study the following particular case of conditional
excess:
$${\bf e}_k^\beta(u) :=\vert \lbrace v\in V(\Gamma):
\partial (u,v)>k  \quad {\rm and }\quad \delta(v)\ge \beta\rbrace\vert .$$
To begin with, firstly we will recall the main properties of the
$k$-alternating polynomials.

The $k$-alternating polynomials, defined and studied in
\cite{alternante} by Fiol, Garriga and Yebra, can be defined as
follows:
 let ${\cal M}= \{\mu_1> \cdots >\mu_b\}$ be a mesh of real numbers. For any
 $k=0,1,...,b-1$ let $P_k$ denote the {\em k-alternating polynomial} associated to  ${\cal M}$.
 That is, the polynomial of  $\mathbb{R}_k[x]$ with norm
 $\| P_k \|_{\infty} =\displaystyle\max_{1 \le i \le b}\{|P_k(\mu_i)| \}$,
such that
 $$ P_k(\mu)= \sup \left\{ P(\mu): P\in \mathbb{R}_k[x], \quad \| P \|_{\infty} \le 1 \right\} $$
 where  $\mu$ is any real number greater than $\mu_1$. We collect here some of its main properties, referring
 the reader to \cite{alternante} for a more detailed study.
  \begin{itemize}
  \item{  For any $k=0,1,...,b-1$ there is a unique $P_k$ which, moreover, is independent of the value of
  $\mu (> \mu_1)$;}
 \item{ $P_k$ has degree k;}
 \item{ $P_0(\mu)=1< P_1(\mu)< \cdots < P_{b-1}(\mu)$;}
 \item{ $P_k$ takes $k+1$  alternating values $\pm 1$ at the mesh
 points;}
 \item{ If $z\in \nu^\perp$ then $  \| P_k({\cal A}(\Gamma))z\| \le \| P_k\|_{\infty} \| z
 \|$ where $\nu=(\sqrt{\delta_1}, \sqrt{\delta_2}, ...,
\sqrt{\delta_n})$, see \cite{degreediameter};}
 \item{ There are
explicit formulae for $P_0(=1),$ $P_1,$ $P_2,$ and $P_{b-1},$
while the other polynomials can be computed by solving a linear
programming problem (for instance by the simplex method).}
 \end{itemize}

\begin{theorem} \label{ThExeceso}
Let $\Gamma=(V,E)$ be a simple and connected graph of size $m$.
Let $u\in V$ and  let $P_k$ be the $k$-alternating polynomial
asso\-ciated to the mesh of the degree-adjacency eigenvalues of
$\Gamma$. Then,
\[
{\bf e}_k^\beta(u)\le
\left\lfloor\frac{2m(2m-\delta(u))}{\beta\left[\delta(u)P_k^2(1)+2m-\delta(u)\right]}\right\rfloor.
\]

\end{theorem}

\begin{proof}
 Let $S=\{v_{i_1}, v_{i_2}, ..., v_{i_s}\}\subset V$ such that $\delta_{i_l}\ge \beta$
$(l=1,...,s)$, and $\partial(S,u) >k$. Let $\sigma=\sum_{l=1}^s
e_{i_l}$, where $e_{i_l}$ denotes the canonical vector
asso\-ciated to the vertex $v_{i_l}$, and let $e$ be the canonical
vector associated to the vertex $u$. From $\partial(S,u) >k
\Rightarrow \langle P_k({\cal A}) \sigma, e \rangle =0$, using the
following decompositions
\begin{equation}\label{desc1}
\sigma=\frac{\langle \sigma,\nu\rangle}{\| \nu\|^2}\nu+
w_s=\frac{\sum_{l=1}^s\sqrt{\delta_{i_l}}}{2m}\nu+ w_s,
\end{equation}
\begin{equation}\label{desc2}
e=\frac{\langle e,\nu\rangle}{\| \nu\|^2}\nu+
w_u=\frac{\sqrt{\delta(u)}}{2m}\nu+ w_u,
\end{equation}
where $\nu=(\sqrt{\delta_1}, \sqrt{\delta_2}, ...,
\sqrt{\delta_n})$  and $w_s,w_u\in\nu^{\perp},$ we obtain
$$P_k(1)\frac{\sqrt{\delta(u)}\sum_{l=1}^s\sqrt{\delta_{i_l}}}{2m}=-\langle
P_k ({\cal A})w_s, w_u \rangle.$$ Hence, by the Cauchy-Schwarz
inequality we have
 $$P_k(1)\frac{\sqrt{\delta(u)}\sum_{l=1}^s\sqrt{\delta_{i_l}}}{2m} \le \| P_k({\cal A}) w_s\| \|
w_u \|.$$ Thus,
\begin{equation}
 P_k(1)\frac{\sqrt{\delta(u)}\sum_{l=1}^s\sqrt{\delta_{i_l}}}{2m} \le \|w_s\| \|w_u\|.\label{ant1}
\end{equation}
 Moreover, the
decompositions (\ref{desc1}) and (\ref{desc2}) lead to
$$  s=\| \sigma\|^2 = \frac{\left(\sum_{l=1}^s\sqrt{\delta_{i_l}}\right)^2}{2m}+ \| w_s\|^2\Rightarrow \| w_s\|=
\sqrt{s-\frac{\left(\sum_{l=1}^s\sqrt{\delta_{i_l}}\right)^2}{2m}}$$
and
$$  1=\| e\|^2 = \frac{\delta(u)}{2m}+ \| w_u\|^2\Rightarrow \| w_u\|=\sqrt{1-\frac{\delta(u)}{2m}}.$$
So, by (\ref{ant1}), we obtain
\begin{equation}\label{final1}
 P_k(1)\sqrt{\delta(u)}\sum_{l=1}^s\sqrt{\delta_{i_l}}\le
\sqrt{(2m-\delta(u))\left(2ms-\left(\sum_{l=1}^s\sqrt{\delta_{i_l}}\right)^2\right)}.
\end{equation}
Therefore,
\begin{equation}\label{final2}
 P_k(1)s\sqrt{\delta(u)\beta}\le
\sqrt{(2m-\delta(u))\left(2ms-s^2\beta\right)}.
\end{equation}
Solving (\ref{final2}) for $s$, and considering that it is an
integer, we obtain the result.
\end{proof}

The above bound is tight for different values of $k$ and $\beta$,
as we can see in the following example. Let $\Gamma$ be the graph
of 5 vertices obtained by joining one vertex of the cycle $C_4$ to
the  vertex of the trivial graph $K_1$. The degree-adjacency
eigenvalues of $\Gamma$ are $\pm 1$, $\pm\frac{\sqrt{6}}{6}$ and
$0$, from which we obtain $P_1(1)=1.84...$ and $P_2(1)=5.899...$.
Hence, the values of the excess ${\bf e}_k^\beta(v)$ are attained
whenever: $\delta(v)=1$, $k=0,1,2$ and $\beta=2,3$; $\delta(v)=2$,
$k=1$ and $\beta=3$; $\delta(v)=2$, $k=2$ and $\beta=1,2,3$;
$\delta(v)=3$, $k=2$ and $\beta=2,3$.

As we can see in  Section \ref{CondWiener}, the above result
becomes an important tool in the study of the conditional Wiener
index.

An analogous upper bound on the standard excess is obtained by
repla\-cing, in above theorem, $\beta$ by the minimum degree
$\delta$. Moreover, in the case of regular graphs,  the above
theorem becomes the following result.

\begin{corollary}
Let $\Gamma$ be a simple and connected graph of order $n$ and let
$P_k$ be the $k$-alternating polynomial associated to the mesh of
the degree-adjacency eigenvalues of $\Gamma$. Then,
\[
{\bf e}_k\le \left\lfloor\frac{n(n-1)}{P_k^2(1)+n-1}\right\rfloor.
\]
\end{corollary}

The above result is analogous to the previous one obtained by the
first author of this paper and Yebra in \cite{exceso}, for
non-necessarily regular graphs,  by using the Laplacian
eigenvalues.

\section{Degree diameter} \label{DegreeDiam}

In this section we study the  problem of finding how far apart can
be two vertices of given degrees in a connected graph. More
precisely, the problem is to find
\[
D^{(\alpha,\beta)}(\Gamma):=\max_{v_i,v_j\in
V}\{\partial(v_i,v_j):\delta_i\ge \alpha, \delta_j\ge \beta \}
\]
 We call this parameter \emph{$(\alpha,\beta)$-degree diameter}.

As in the case of the standard diameter, the study of this
parameter is of interest in the design of interconnection networks
when we need to minimize the communication delays between  two
nodes of given degrees.

In this section we obtain a tight bound on the
$(\alpha,\beta)$-degree diameter by using the $k$-alternating
polynomials  on the mesh of eigenvalues of the degree-adjacency
matrix.

\begin{theorem} \label{theorem}
Let $\Gamma=(V,E)$ be a simple and connected graph of size $m$.
Let $P_k$ be the $k$-alternating polynomial associated to the mesh
of the degree-adjacency eigenvalues of $\Gamma$.  Then,
\begin{equation}\label{cotaDegreeDiam}
P_k(1)
>\sqrt{\left(\frac{2m}{\alpha}-1\right)\left(\frac{2m}{\beta}-1\right)}
  \Rightarrow D^{(\alpha,\beta)}(\Gamma) \le k.
\end{equation}
\end{theorem}

\begin{proof}
 Let $e_i$ and $e_j$ be the canonical vectors of $\mathbb{R}^n$
associated to the vertices $v_i$ and $v_j$. Using the following
decomposition
\begin{equation}\label{desc}
e_i=\frac{\langle e_i,\nu\rangle}{\| \nu\|^2}\nu+
u=\frac{\sqrt{\delta_i}}{2m}\nu+ u, \quad e_j=\frac{\langle
e_j,\nu\rangle}{\| \nu\|^2}\nu+ w=\frac{\sqrt{\delta_j}}{2m}\nu+
w,
\end{equation}
where $\nu=(\sqrt{\delta_1}, \sqrt{\delta_2}, ...,
\sqrt{\delta_n})$ and $u,w\in\nu^{\perp},$ we obtain
\begin{align*}
\partial(v_i,v_j) >k &\Rightarrow \left(P_k\left( {\cal A}
\right)\right)_{ij}=0\\
 &\Rightarrow \langle P_k({\cal A}) e_i, e_j \rangle =0 \\
  &\Rightarrow P_k(1)\frac{\sqrt{\delta_i\delta_j}}{2m}+  \langle P_k({\cal A}) u, w
  \rangle=0\\
    &\Rightarrow P_k(1)\frac{\sqrt{\delta_i\delta_j}}{2m}=-\langle P_k ({\cal A})u, w \rangle.
\end{align*}
Then, by the Cauchy-Schwarz inequality we have

 \begin{align}
\partial(v_i,v_j) >k &\Rightarrow P_k(1)\frac{\sqrt{\delta_i\delta_j}}{2m} \le \| P_k({\cal A}) u\| \| w \|\\
                    &\Rightarrow  P_k(1)\frac{\sqrt{\delta_i\delta_j}}{2m} \le
                        \|P_k\|_{\infty} \|u\| \|w\|.\label{ant}
\end{align}
Moreover, the decomposition (\ref{desc}) leads to
$$  1=\| e_i\|^2 = \frac{\delta_i}{2m}+ \| u\|^2\Rightarrow \| u\|=\sqrt{1-\frac{\delta_i}{2m}}$$
and
$$  1=\| e_j\|^2 = \frac{\delta_j}{2m}+ \| w\|^2\Rightarrow \| w\|=\sqrt{1-\frac{\delta_j}{2m}}$$
So, by (\ref{ant}), we obtain
\begin{equation}\label{final}
\partial(v_i,v_j) >k \Rightarrow  P_k(1)\sqrt{\delta_i\delta_j}\le
\sqrt{(2m-\delta_i)(2m-\delta_j)}  .
\end{equation}
The converse of (\ref{final}) leads to
\begin{equation}\label{fund}
P_k(1) >\sqrt{\frac{(2m-\delta_i)(2m-\delta_j)}{\delta_i
\delta_j}}
 \Rightarrow  \partial(v_i,v_j) \le k.
\end{equation}
The result follows from (\ref{fund}).
\end{proof}

\begin{figure}[h]
\begin{center}
\caption{  }\label{ejemplo}
\includegraphics[height=3cm]{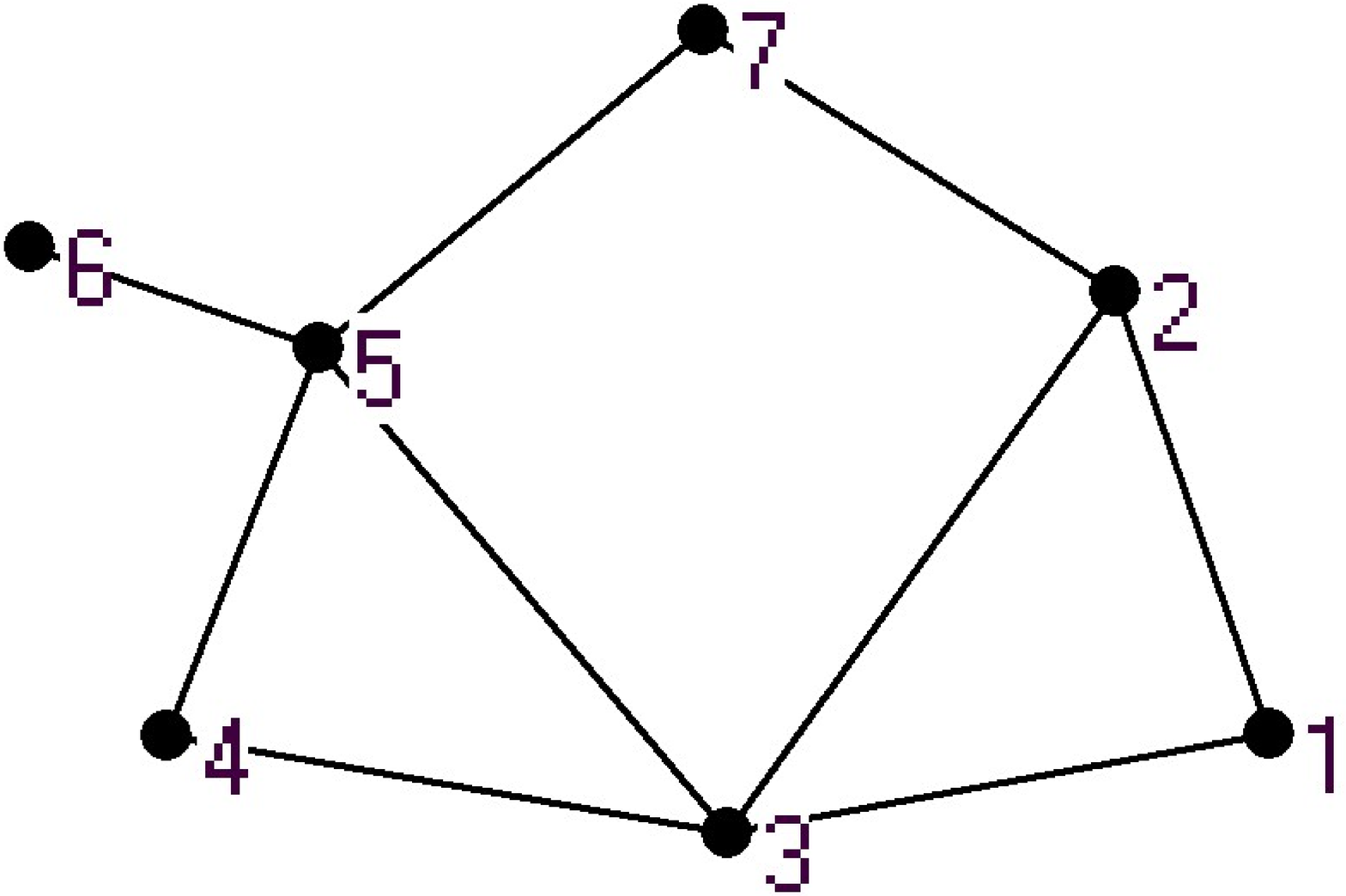}
\end{center}
\end{figure}

As we can see in the following example, the above bound is
attained for several values of $\alpha$ and $\beta$.
 The graph of Figure \ref{ejemplo} has degree-adjacency eigenvalues
\begin{center}
$\left\{1,
\frac{-3+\sqrt{249}}{24},\frac{1}{4},0,-\frac{1}{2},-\frac{1}{2},\frac{-3-\sqrt{249}}{24}\right\}$
\end{center}
from which we obtain
\begin{center}
 $P_1(1)=1.7$, $P_2(1)=5$, $P_3(1)=15.2$ and $P_4(1)=58$.
\end{center}
Thus, the following bounds are attained:
\begin{center}
 $D^{(1,2)}(\Gamma)\le 3$, $D^{(3,4)}(\Gamma) \le 2$ and $D^{(4,4)}(\Gamma) \le 1$.
\end{center}

As particular cases of above theorem we derive the following
results in which the expression (\ref{cotaDegreeDiam}) is
simplified.

\begin{corollary}
Let $\Gamma=(V,E)$ be a simple and connected graph of order $n$
and size $m$. Let $P_k$ be the $k$-alternating polynomial
associated to the mesh of the degree-adjacency eigenvalues of
$\Gamma$. Then,

\begin{equation} \label{betabeta}
 P_k(1)>\frac{2m}{\beta}-1 \Rightarrow D^{(\beta,\beta)}(\Gamma) \le k.
 \end{equation}

The standard diameter is bounded by
\begin{equation}
P_k(1)>\frac{2m}{\delta}-1 \Rightarrow D(\Gamma) \le k.
\end{equation}

If $\Gamma$ is regular, the standard diameter is bounded by
\begin{equation} \label{casoReg}
P_k(1)>n-1 \Rightarrow D(\Gamma) \le k.
\end{equation}
\end{corollary}

 As we can see in next section, the bound (\ref{betabeta}) becomes an important tool in the study
of the conditional Wiener index. Moreover, the bound
(\ref{casoReg}) is an analogous result to the previous one given
by Fiol, Garriga and Yebra in \cite{alternante} by using the
standard adjacency matrix. The reader is referred to
\cite{degreediameter} for a more general study on the conditional
diameter.

\section{Conditional Wiener index}\label{CondWiener}

The {\em Wiener index} $W(\Gamma)$ of a graph $\Gamma$ with vertex
set $\{v_1,v_2,...,v_n\}$  defined as the sum of distances between
all pairs of vertices of $\Gamma$,
 $$W(\Gamma):= \frac{1}{2}\sum_{i=1,j=1}^n\partial(v_i,v_j),$$
 is the first mathematical invariant reflecting the
topological structure of a molecular graph.

This topological index  has been extensively studied, for
instance, a comprehensive survey on the direct calculation,
applications and the relation of the Wiener index of trees with
other parameters of graphs can be found in \cite{wiennerofTrees}.
Moreover, a list of 120 references of the main works on the Wiener
index of graphs can be found in the referred survey.

Alternatively, the Wiener index can be defined as
$$W(\Gamma)= \frac{1}{2}\sum_{v\in V(\Gamma)}S(v),$$
where $S(v)$ denotes the {\em distance of the vertex} $v$:
$$S(v):=\sum_{u\in V(\Gamma)}\partial(u,v).$$
We define the {\em conditional Wiener index}
$$W_\wp(\Gamma):=\frac{1}{2}\sum_{v\in \wp} S_{\wp}(v),$$
where $\wp$ is a property and $v\in \wp$ means that the vertex $v$
satisfies the property $\wp$, and
$$S_\wp(v):=\sum_{u\in \wp}\partial(u,v)$$
is the {\em conditional distance} of $v$. In particular, if $\wp$
requires that $\delta(v)\ge \beta$, the conditional Wiener index
will be denoted by $W_\beta(\Gamma)$, moreover, the  conditional
distance of $v$ will be denoted by $S_\beta (v)$. Clearly, if
$\beta$ is the minimum degree of $\Gamma$, then $W_\beta(\Gamma)$
and the standard Wiener index coincides.

\begin{lemma}\label{pdemostrar}
The conditional  Wiener index  of a graph  $\Gamma$,
$W_\beta(\Gamma)$, satisfies
\[
 W_\beta(\Gamma) = \frac{1}{2} \sum_{\delta(v)\ge \beta} \sum_{k= 0}^{D^{(\beta,\beta)}(\Gamma)-1} {\bf e}_k^\beta(v)
\]
\end{lemma}

\begin{proof}
 For each vertex $v\in V(\Gamma)$ of degree $\delta(v)>\beta$  we have
$$S_\beta(v)=\sum_{k=1}^{D^{(\beta,\beta)}(\Gamma)}k({\bf e}^\beta_{k-1}(v)- {\bf e}^\beta_{k}(v)).$$
Moreover, by a simple calculation we have
\begin{equation} \label{distvertex}
S_\beta(v)=\sum_{k=0}^{D^{(\beta,\beta)}(\Gamma)-1}{\bf
e}^\beta_{k}(v).
\end{equation}
Hence, by (\ref{distvertex}) we obtain the result.
\end{proof}

 Therefore, it follows from Lemma \ref{pdemostrar} that bounds on ${\bf
e}^\beta_{k}$ lead to bounds on the conditional Wiener index
$W_\beta$.

\begin{theorem}
Let $\Gamma=(V,E)$ be a simple and connected graph of size $m$.
Let $P_k$ be the $k$-alternating polynomial associated to the mesh
of the degree-adjacency eigenvalues of $\Gamma$ and let
$x=\left|\{v\in V(\Gamma): \quad \delta(v)\ge \beta\}\right|$. If
$P_k(1)> \frac{2m}{\beta}-1$, then
$$W_\beta(\Gamma)\le \frac{x}{2}\displaystyle\sum_{l=0}^{k-1}\left\lfloor\frac{2m(2m-\beta)}
{\beta\left(\beta P_l^2(1)+2m-\beta\right)}\right\rfloor.$$
\end{theorem}

\begin{proof}
By Lemma \ref{pdemostrar} and Theorem \ref{ThExeceso} we have
\begin{equation}
W_\beta(\Gamma)\le
\frac{x}{2}\displaystyle\sum_{k=0}^{D^{(\beta,\beta)(\Gamma)-1}}\left\lfloor\frac{2m(2m-\beta)}
{\beta\left(\beta P_k^2(1)+2m-\beta\right)}\right\rfloor.
\end{equation}
Therefore, by (\ref{betabeta}) we conclude the proof.
\end{proof}

An analogous upper bound on the standard Wiener index is obtained
by replacing, in above theorem, $\beta$ by $\delta$, and $x$ by
$n$. Moreover, in the case of regular graphs, the above theorem
becomes the following result.

\begin{corollary}
Let $\Gamma$ be a simple and connected $\delta$-regular graph of
order $n$. Let $P_k$ be the $k$-alternating polynomial associated
to the mesh of the degree-adjacency eigenvalues of $\Gamma$. If
$P_k(1)> n-1$, then
$$W(\Gamma)\le \frac{n}{2}\displaystyle\sum_{l=0}^{k-1}\left\lfloor\frac{n(n-1)}
{P_l^2(1)+n-1}\right\rfloor.$$
\end{corollary}
The reader is referred to \cite{WienerHypergraph} for a more
general study on the Wiener index of hypergraphs.


\begin{thebibliography}{99}







\bibitem{Biggs}
N. Biggs, {\it Algebraic graph theory}, Cambridge University
Press, 1993.

\bibitem{Bollobas} B. Bollob\'{a}s and P. Erd\"{o}s, Graphs of
extremal weights, {\it Ars Combinatoria}  {\bf 50} (1998),
225-233.



\bibitem{cvetkovic}
D. M. Cvetkovi\'{c}, M. Doob and H. Sachs, {\it Spectra of
graphs}, Academic Press Inc., New York, 1979.


\bibitem{Zagreb}
K. C. Das and I. Gutman, Some properties of the second Zagreb
index, {\it MATCH Commun.  Math.  Comput. Chem.} {\bf 52} (2004),
103-112.

\bibitem{wiennerofTrees}
A. A. Dobrynin,   R. Entringer   and I. Gutman,   Wiener Index of
Trees: Theory and Applications, {\it Acta Applicandae
Mathematicae} {\bf 66}  (2001), 211-249.


\bibitem{alternante}
  M.A. Fiol, E. Garriga and J.L.A. Yebra, On a class of polynomials and its relation with the
spectra and diameters of graphs, {\it  J. Combin. Theory Ser. B}
{\bf  67} (1996), 48-61.







\bibitem{from-local}
M.A. Fiol and E. Garriga,  From local adjacency polynomials to
locally pseudo-distance-regular-graphs. {\it   J. Combin. Theory
Ser. B}. {\bf 71} (1997),  162-183.




\bibitem{Godsil}
C. Godsil and G. Royle, {\it Algebraic graph theory},
Springer-Verlang New York, Inc. 2001.

\bibitem{gutman}
I. Gutman and M. Lepovi\'{c}, Choosing the exponent in the
definition of the connectivity, {\it J. Serb. Chem. Soc.} {\bf 66}
(2001), 605-611.

\bibitem{Randic} M. Randi\'{c}, On the characterization of
molecular branching. {\it J. Amer. Chem. Soc.} {\bf 97} (1975),
6609-6615.



\bibitem{exceso} J. A. Rodr\'{\i}guez
and  J.L.A. Yebra, Laplacian eigenvalues and the excess of a graph
{\it Ars Combinatoria.} {\bf 64} (2002),  249-258.


\bibitem{eama}
J. A. Rodr\'{\i}guez and J. L. A. Yebra, Cotas espectrales del
exceso y la excen\-tricidad de un conjunto de v\'ertices  de un
grafo.  {\it EAMA-97} (1997), 383-390.

\bibitem{degreediameter}  J. A. Rodr\'{\i}guez, The $(\alpha,\beta, s,t)$-diameter of  graphs: A particular case
of conditional diameter. \textit{Discrete Applied Mathematics }\textbf{154} (14) (2006) 2024--2031.


\bibitem{Randic1}
J. A. Rodr\'{\i}guez, A spectral approach to the Randi\'{c} index. \textit{Linear Algebra and its Applications} \textbf{400} (2005) 339--344.

\bibitem{WienerHypergraph}  J. A. Rodr\'{\i}guez, On the Wiener index and the eccentric distance sum of
hypergraphs,  {\it MATCH Commun.  Math. Comput. Chem.} {\bf 54}
(1) (2005)  209--220.

\end{thebibliography}
\end{document}